\documentclass[oneside,11pt]{amsart}
\usepackage{amsmath, amsfonts,amsthm,times,graphics}
\usepackage[active]{srcltx}

\newtheorem{theorem}{Theorem}[section]
\newtheorem{lemma}[theorem]{Lemma}

\newtheorem{proposition}[theorem]{Proposition}

\newtheorem{remark}[theorem]{Remark}
\newtheorem{conjecture}[theorem]{Conjecture}
\begin{document}
\title{Norm of the Bergman projection}
\author{David Kalaj}
\address{
Faculty of Natural Sciences and Mathematics, University of
Montenegro, Cetinjski put b.b. 81000 Podgorica, Montenegro}
\email{davidk@ac.me}

\author{Marijan Markovi\'c}
\address{
Faculty of Natural Sciences and Mathematics, University of
Montenegro, Cetinjski put b.b. 81000 Podgorica, Montenegro}
\email{marijanmmarkovic@gmail.com}


\begin{abstract}
This paper deals with the the norm of the weighted Bergman
projection operator $P_\alpha:L^\infty(\mathbb{B})\to \mathcal{B}$
where $\alpha>-1$ and $\mathcal{B}$ is the Bloch space of the unit
ball $\mathbb{B}$ of the complex space $\mathbb{C}^n$. We consider
two Bloch norms, the standard Bloch norm and invariant norm w.r.t.
automorphisms of the unit ball. Our work contains as a special case
the main result of the recent paper \cite{Perala}.
\end{abstract}
\maketitle

\section{Introduction and preliminaries}
Introduce the notation which will be used in this paper. We follow
the Rudin monograph \cite{Rudin-Book-Ball}. Throughout the paper $n$
is an integer bigger or equal to $1$. Let $\left<\cdot,\cdot\right>$
stands for the inner product in the complex $n$-dimensional space
$\mathbb C^n$ given by
$$\left<z,w\right>=z_1\overline{w}_1+\dots+z_n\overline{w}_n,\quad z,\ w\in\mathbb C^n,$$
where $z=(z_1,\dots,z_n)$ and $w=(w_1,\dots,w_n)$ are coordinate
representation in the standard base $\{e_1,\dots,e_n\}$ of
$\mathbb{C}^n$. The inner product induces the Euclidean norm
$$|z|=\left<z,z\right>^{1/2},\quad z\in\mathbb C^n.$$
Denote by ${\mathbb B}$ the unit ball $\{z\in\mathbb C^n:|z|<1\}$
and let ${\mathbb S}=\partial {\mathbb B}$ be its boundary.

We let $dv$ be the volume measure in $\mathbb {C}^n$, normalized so
that $v({\mathbb B})=1$. We will also consider a class of weighted
volume measures on ${\mathbb B}$. When $\alpha>-1$, we define a
finite measure $dv_\alpha$ on ${\mathbb B}$ by
$$dv_\alpha(z)=c_\alpha (1-|z|^2)^\alpha dv(z),$$
where $c_\alpha$ is a normalizing constant so that
$v_\alpha({\mathbb B})=1$. Using polar coordinates, one can easily
calculate that \begin{equation}\label{cal}c_\alpha=
\binom{n+\alpha}{n}.\end{equation}

It is known that the biholomorphic mappings of ${\mathbb B}$ onto
itself have the following form
$$\varphi_a(\omega)=
\frac
{a-\frac{\left<\omega,a\right>}{|a|^2}a-(1-|a|^2)^{1/2}(\omega-\frac{\left<\omega,a\right>}{|a|^2}a)}
{1-\left<\omega,a\right>}, \quad \text{for}\quad a\in {\mathbb B},$$
up to unitary transformations; for $a=0$, we set $\varphi_a=-
\mathbf{Id}_{\mathbb B}$. In the case $n=1$ this is simply the
equality $\varphi_a(\omega)= (a-\omega)/(1-\overline{a}\omega)$.
Traditionally, these mappings are also called biholomorphic
automorphisms. By $\mathbf{Aut}({\mathbb B}) =
\{U\circ\varphi_a:a\in {\mathbb B},\ U\in \mathcal{U}\}$, where
$\mathcal{U}$ is the group of all unitary transformations of the
space $\mathbb {C}^n$, is denoted the group of all biholomorphic
automorphisms of the unit ball. One often calls
$\mathbf{Aut}({\mathbb B})$ the group of M\"obius transformations of
$\mathbb B$.

Observe that $\varphi_a(0)=a$. Since $\varphi_a$ is involutive, i.e.
$\varphi_a\circ\varphi_a={\mathbf{Id}}_{\mathbb B}$, we also have
$\varphi_a(a)=0$.

Viewing $\mathbb {C}^n$ as $\mathbb {R}^{2n}$, the real Jacobian of
$\varphi_a$ is given by
$$(J_\mathbb R\varphi_a)(\omega)=\left(\frac{1-|a|^2}{\left|1-\left<\omega,a\right>\right|^2}\right)^{n+1},\quad \omega\in {\mathbb B}.$$
Two identities
\begin{equation}\label{Relation1}
1-|\varphi_a(\omega)|^2=\frac{(1-|a|^2)(1-|\omega|^2)}{|1-\left<\omega,a\right>|^2}
\end{equation}
and
\begin{equation}\label{Relation2}
(1-\left<\omega,a\right>)(1-\left<\varphi_a(\omega),a\right>)=1-|a|^2,
\end{equation}
for all $a,\ \omega\in B$, will also be useful. By using
\eqref{Relation1} we obtain the next relation
\[\begin{split}
dv_\alpha(\varphi_a(\omega))& =(1-|\varphi_a(\omega)|^2)^\alpha
(J_\mathbb R\varphi_a)(\omega) dv(\omega)
\\&=\left(\frac{(1-|\omega|^2)(1-|a|^2)}{|1-\left<\omega,a\right>|^2}\right)^\alpha
\left(\frac{1-|a|^2}{\left|1-\left<\omega,a\right>\right|^2}\right)^{n+1}
dv(\omega)
\\&=\left(\frac{\left(1-|a|^2\right)}{\left|1-\left<\omega,a\right>\right|^2}\right)^{n+1+\alpha} dv_\alpha(\omega).
\end{split}\]
For a holomorphic function $f$ with $\nabla f$ we denote the complex
gradient
$$\nabla f(z)=\left(\frac{\partial f}{\partial z_1}(z),\dots,\frac{\partial f}{\partial z_n}(z)\right).$$

The Bloch space $\mathcal{B}$ contains all functions $f$ holomorphic
in ${\mathbb B}$ for which the semi-norm
$$\|f\|_{\beta}:=\sup_{z\in {\mathbb B}}(1-|z|^2)\left|\nabla f(z)\right|$$
is finite. One can obtain a true norm by adding $|f(0)|$, more
precisely in the following way
$$\|f\|_\mathcal{B}=|f(0)|+\|f\|_\beta,\quad f\in\mathcal{B}.$$
It is well known that $\mathcal{B}$ is a Banach space with the above
norm. The standard reference for Bloch space of the unit disc is
\cite{Anderson}. For the high dimension case we refer to
\cite{timo-BLMS}, \cite{timo-Crelle} and \cite{Zhu-Book-Spaces}.

Let $L^p$ stands for Lebesgue space of all measurable functions in
$\mathbb B$ which modulus is integrable in $\mathbb B$ with exponent
$p$ when $1\le p<\infty$ and for $p=\infty$ the space of essentially
bounded measurable functions in the unit ball. The Bergman
projection operator $P_\alpha$ for $\alpha>-1$ is defined by
$$P_\alpha g(z)=\int_{\mathbb B}\mathcal{K}_\alpha(z,w)g(w)dv_\alpha(w),\quad g\in L^p(\mathbb B),$$
where
$$\mathcal{K}_\alpha(z,w)=\frac 1{(1-\left<z,w\right>)^{n+1+\alpha}},\quad z,\ w\in {\mathbb B}$$
is the weighted Bergman kernel. Bergman type projections are central
operators when dealing with questions related to analytic function
spaces. One often wants to prove that Bergman projections are
bounded and the exact operator norm of the operator is difficult to
obtain. By Forelli--Rudin theorem, $P_\alpha$ is bounded if and only
if $\alpha>1/p-1$. Thus, it is not bounded as operator $L^1\to A^1$
and it is known that it is not bounded as $L^\infty\to H^\infty$. On
the other hand, for $n=1$ it is well known that the Bergman
projection $P_\alpha: L^{\infty}({\mathbb B})\to\mathcal{B}$ is
bounded and onto, see \cite{Zhu-Book-Operator}. For $n>1$ the
operator $P_\alpha: L^{\infty}({\mathbb B})\to \mathcal{B}$ is
surjective what can be seen from
\cite[Theorem~3.4]{Zhu-Book-Spaces}) in the Zhu book.

The $\beta-$norm and $\mathcal{B}-$norm of the Bergman projection
$P_\alpha:L^{\infty}({\mathbb B})\to\mathcal{B}$ are
$$\|P_\alpha\|_\beta=\sup_{\|g\|_\infty\le 1}\|P_\alpha g\|_\beta,
\quad\text{and}\quad \|P_\alpha\|_\mathcal{B}=\sup_{\|g\|_\infty\le
1}\|P_\alpha g\|_\mathcal{B}.$$ There are several equivalent ways to
introduce the Bloch spaces in the ball ${\mathbb
B}\subseteq\mathbb{C}^n$. The previous one is natural and
straightforward but the norm defined in that way is not invariant
with respect to the group $\mathbf{Aut}({\mathbb B})$. The following
Bloch norm has this property.

We define the invariant gradient $|\tilde \nabla f(z)|$ where
$$\tilde\nabla f(z)=\nabla (f\circ \varphi_z)(0),$$ where $\varphi_z$
is an automorphisms of the unit ball such that $\varphi_z(0)=z$.
This norm is invariant w.r.t. automorphisms of the unit ball. Namely
$$|\tilde\nabla (f\circ \varphi)|=|(\tilde\nabla f)\circ \varphi|$$
for $\varphi\in \mathbf{Aut}({\mathbb B})$. Then the Bloch space
$\mathcal{B}$ contains all holomorphic functions $f$ in the ball
${\mathbb B}$ for which
$$\|f\|_\mathcal{\tilde\beta}:=\sup_{z\in {\mathbb B}}
|\tilde\nabla f(z)|<\infty$$ (cf.
\cite[Theorem~3.4]{Zhu-Book-Spaces} or \cite{timo-BLMS}). For $n=1$
we have $|\tilde\nabla f(z)|=(1-|z|^2)\left|\nabla f(z)\right|$, but
for $n>1$ this is not true. Notice that
$\|\cdot\|_\mathcal{\tilde\beta}$ is also a semi-norm. One can
obtain a norm in the following way
$$\|f\|_{\tilde{\mathcal{B}}}=|f(0)|+\|f\|_{\tilde\beta},\quad f\in\mathcal{B}.$$
The  $\tilde\beta$-norm ($\tilde{\mathcal{B}}$-norm) of the Bergman
projection is
$$\|P_\alpha\|_{\tilde{\beta}}=\sup_{\|f\|\le 1}\|P_\alpha
g\|_{\tilde{\beta}}, \quad
(\|P_\alpha\|_{\tilde{\mathcal{B}}}=\sup_{\|g\|\le 1}\|P_\alpha
g\|_{\tilde{\mathcal{B}}}).$$ From the proof of
\cite[Theorem~3.4]{Zhu-Book-Spaces}) we find out that
$$\| P_\alpha g\|_{\tilde{\beta}}\le C\|g\|_\infty,$$
where $C$ is a positive constant. The later implies that $P_\alpha$
is a bounded operator since
$$\|P_\alpha\|_{\tilde{\mathcal{B}}}\le 1+\|P_\alpha\|_{\tilde\beta}.$$

Before stating the main results let us prove the following simple
lemma.

\begin{lemma}\label{mele}
For $\alpha>-1$ we have
\begin{equation}\label{prima}\|P_\alpha\|_{\mathcal{B}}\le
1+\|P_\alpha\|_{\beta}\end{equation} and
\begin{equation}\label{seconda}\|P_\alpha\|_{\tilde{\mathcal{B}}}\le
1+\|P_\alpha\|_{\tilde\beta}.\end{equation}
\end{lemma}
\begin{proof}
Since $$|P_\alpha g(0)|=\left|\int_{\mathbb B}
g(w)dv_\alpha(w)\right|\le \|g\|_\infty$$ it follows that
$$\|P_\alpha g\|_\mathcal{B}=|P_\alpha g(0)|+\|P_\alpha g\|_\beta\le
\|g\|_\infty + \|P_\alpha\|_\beta \|g\|_\infty.$$ This implies
\eqref{prima}. The relation \eqref{seconda} can be proved similarly.
\end{proof}
In this paper we find the exact norm of $P_\alpha$ w.r.t.
$\beta-$Bloch (semi) norm, and estimate the $\tilde \beta-$Bloch
(semi) norm. It is the content of our Theorem~\ref{Main-Theorem}
which generalizes the result from the recent paper \cite{Perala} in
two directions. Let
$$C_{\alpha,n}:= \frac{\Gamma(2 + n+\alpha)
}{\Gamma^2((2+n+\alpha)/2)}.$$ In this paper we prove the following
two theorems

\begin{theorem}\label{Main-Theorem}
For the $\beta-$(semi) norm of the Bergman projection $P_\alpha$ we
have
$$\|P_\alpha\|_\beta=C_{\alpha,n}$$  where
$\Gamma$ is Euler's Gamma function.
\end{theorem}
In order to formulate the next theorem, assume that $n>1$ and define
\begin{equation}\label{ell}\ell(t)=(n+1+\alpha)\int_{\mathbb{B}}\frac{|(1-w_1)\cos
t+w_2\sin t|}{|w_1-1|^{n+1+\alpha}}dv_\alpha(w).\end{equation}
\begin{theorem}\label{Mait}
For $\alpha>-1$ we have
\begin{equation}\label{elles}\ell(\pi/2)=\frac{\pi}{2} \ell (0)=\frac{\pi}2
C_{\alpha,n}.\end{equation}
For the $\tilde \beta-$ (semi) norm of
the Bergman projection $P_\alpha$ we have
\begin{equation}\label{esti1}\|P_\alpha\|_{\tilde
\beta}=\tilde{C}_{\alpha,n}=\max_{0\le
t\le\pi/2}\ell(t)\end{equation} and
\begin{equation}\label{esti2} \frac{\pi}{2}C_{\alpha,n}\le \|P_\alpha\|_{\tilde \beta}\le
\frac{\sqrt{\pi^2+4}}{2} C_{\alpha,n}.\end{equation}
\end{theorem}
\begin{remark}
For $\alpha=0$ we put $P=P_\alpha$ and we have $$\|P\|_\beta= \frac
{(n+1)!} { \Gamma^2(1+n/2)}.$$ Moreover, for $n=1$ we obtain
$$\|P\|_\beta=\frac 8\pi$$ which
presents the main result in \cite{Perala}. As an immediate corollary
of Theorem~\ref{Main-Theorem}, Theorem~\ref{Mait} and
Lemma~\ref{mele} we have the next norm estimates of the Bergman
projection:
$$C_{\alpha,n}\le \|P_\alpha\|_\mathcal{B}\le 1+C_{\alpha,n}$$
and
$$\frac{\pi}{2}C_{\alpha,n}\le \|P_\alpha\|_{\tilde{\mathcal{B}}}\le 1+\frac{\sqrt{\pi^2+4}}{2}C_{\alpha,n}.$$
\end{remark}
\begin{conjecture}\label{cok}
In connection with Theorem~\ref{Mait}, we conjecture that
$$\tilde{C}_{\alpha,n}=\frac{\pi}{2}C_{\alpha,n}.$$ See
\textbf{Appendix}
below for an approach that can be of interest.
\end{conjecture}

\section{Proof of the Theorem~\ref{Main-Theorem}}
What we have to find is
$$\|P_\alpha\|_\beta=\sup \{ (1-|z|^2)\left|\nabla_z(P_\alpha g)(z)\right|  : |z|<1,\ \|g\|_\infty\le 1 \}.$$
A straightforward calculation yields
\begin{equation}\label{nabla}
\nabla_z\mathcal{K}_\alpha(z,w)=\frac{(1+n+\alpha)\overline{w}}{(1-\left<z,w\right>)^{n+2+\alpha}},\quad
z,\ w\in {\mathbb B},
\end{equation}
and this implies the formula
$$\nabla_z(P_\alpha g)(z)=\int_{\mathbb B} \nabla_z \mathcal{K}_\alpha(z,w) g(w)dv_\alpha(w),\quad z\in {\mathbb B}.$$
First of all, for a fixed $z\in {\mathbb B}$ and for
$\|g\|_\infty\le 1$ we have the following estimates
\[\begin{split}
\left|\nabla( P_\alpha g)(z)\right|& =\max_{\zeta\in {\mathbb S}}
\left|\left<\nabla P_\alpha g(z),\zeta\right>\right|
\\&=\max_{\zeta\in {\mathbb S}}
\left|\left<\int_{\mathbb B} \nabla_z
\mathcal{K}_\alpha(z,w)g(w)dv_\alpha(w),\zeta\right>\right|
\\&=\max_{\zeta\in {\mathbb S}}
\left|\int_{\mathbb B}
\left<\nabla_z\mathcal{K}_\alpha(z,w)g(w),\zeta\right>dv_\alpha(w)\right|
\\&\le \max_{\zeta\in {\mathbb S}}
\int_{\mathbb B}
\left|\left<\nabla_z\mathcal{K}_\alpha(z,w)g(w),\zeta\right>\right|dv_\alpha(w)
\\&=\max_{\zeta\in {\mathbb S}}
\int_{\mathbb
B}\left|\left<\frac{(1+n+\alpha)\overline{w}}{(1-\left<z,w\right>)^{n+2+\alpha}},\zeta\right>\right||g(w)|dv_\alpha(w)
\\&\le\max_{\zeta\in {\mathbb S}}
\int_{\mathbb
B}\frac{(1+n+\alpha)\left|\left<\overline{w},\zeta\right>\right|}{\left|1-\left<z,w\right>\right|^{n+2+\alpha}}dv_\alpha(w).
\end{split}\]
Denote
$$F_\zeta(z)
=(1+n+\alpha)\int_{\mathbb B} \frac
{(1-|z|^2)\left|\left<w,\overline{\zeta}\right>\right|}{\left|1-\left<z,w\right>\right|^{n+2+\alpha}}
dv_\alpha(w).$$ The statement of the Theorem~\ref{Main-Theorem} will
follow directly from the following equalities
$$\|P_\alpha\|_\beta=\sup\{F_\zeta(z):z\in {\mathbb B},\ \zeta\in {\mathbb S}\}=C_{\alpha},$$
which will be proved through the following two lemmas:

\begin{lemma}\label{le1}
For every $\alpha>-1$ we have
$$\sup\{F_\zeta(z):z\in {\mathbb B},\ \zeta\in {\mathbb S}\}\le C_{\alpha,n}.$$
\end{lemma}

\begin{lemma}\label{le2}
For every $\alpha>-1$ there exists a sequence $g_k$ of functions
$\|g_k\|_\infty=1$ and a sequence of vectors $z_k\in {\mathbb B},\
k\ge1 $ such that
$$\lim_{k\to \infty}(1-|z_k|^2)|\nabla( P_\alpha g_k)(z_k)|=C_{\alpha,n}.$$
\end{lemma}
In order to give proofs of the previous lemmas we need
\cite[Proposition~1.4.10]{Rudin-Book-Ball} and some its corollaries
collected in the following proposition.
\begin{proposition}\label{rudin}
a) For $z\in {\mathbb B}$, $c$ real, $t>-1$ define
$$J_{c,t}(z)=\int_{\mathbb B}\frac{(1-|w|^2)^t}{|1-\left<z,w\right>|^{n+1+t+c}}dv(w).$$
When $c<0$, then $J_{c,t}$ is bounded in ${\mathbb B}$. Moreover,
\begin{equation}\label{rud}J_{c,t}(z)=\frac{\Gamma(n+1)\Gamma(1+t)}{\Gamma(\lambda_1)^2}
\sum_{k=0}^\infty\frac{\Gamma^2(k+\lambda_1)|z|^{2k}}{\Gamma(k+1)\Gamma(n+1+t+k)},\end{equation}
where $\lambda_1=\frac{n+1+t+c}{2}$.
\\
b) Further we can write $J_{c,t}$ in the closed form as
\begin{equation}\label{close} J_{c,t}(z)=\frac{\Gamma(1 + n)
\Gamma(1 + t)\ F[{\lambda_1, \lambda_1}, {1 + n + t},
  |z|^2]}{\Gamma(1 + n + t)},\end{equation} where $F$ is the Gauss hypergeometric function. In particular
  \begin{equation}\label{spec}J_{c,t}\left(\frac{z}{|z|}\right)=\frac{\Gamma(1 + n) \Gamma(1 + t) \Gamma(-c)}{\Gamma^2(
 1/2 (1 - c + n + t))}.
\end{equation}
\end{proposition}
\begin{proof}[Proof of Proposition~\ref{rudin}]
The first part of proposition coincides with the first part of
\cite[Proposition~1.4.10]{Rudin-Book-Ball} together with its proof.
In order to prove the part b) we recall the classical definition of
the Gauss hypergeometric function:
\begin{equation}\label{defh}{F}(a, b, c, z) = 1 +
\sum_{n=1}^\infty\frac{(a)_n(b)_n}{ (c)_n n!} z^n,\end{equation}
where $(d)_n=d(d+1)\cdots (d+n-1)$ is the Pochhammer symbol. The
series converges at least for complex $z \in
\mathbb{U}:=\{z:|z|<1\}\subset \mathbb C$ and for $z\in
\mathbb{T}:=\{z:|z|=1\}$, if $c>a+b$. For $\Re(c)>\Re (b)>0$ we have
the following well-known formula
\begin{equation}\label{for}F(a,b,c,z)=\frac{\Gamma(c)}{\Gamma(b)\Gamma(c-b)}\int_0^1\frac{t^{b-1}(1-t)^{c-b-1}}{(1-tz)^a}dt.\end{equation}
In particular the Gauss theorem states that
\begin{equation}\label{fora}F(a,b,c,1)=\frac{\Gamma(c)\Gamma(c-a-b)}{\Gamma(b)\Gamma(c-b)},\ \ \Re (c)>\Re(a+b).\end{equation}
In order to derive \eqref{close} from \eqref{rud}, we use the
formula $\Gamma(x+1)=x\Gamma(x)$ and obtain
\begin{equation}\label{shtu}\Gamma(k+\lambda_1)=(\lambda_1)_k \Gamma(\lambda_1)\end{equation} and
\begin{equation}\label{die}\Gamma(n+1+t+k)=(n+1+t)_k\Gamma(n+1+t).\end{equation}
From \eqref{defh}, \eqref{shtu} and \eqref{die}, by taking
$a=b=\lambda_1$ and $c=1+n+t$,  we derive \eqref{close}. The formula
\eqref{spec} follows from \eqref{fora} and observing that
$c>a+b=1+n+t+c$.

\end{proof}
Also we need the Vitali theorem, and include its formulation (cf.
\cite[Theorem 26.C]{halmos}).

\begin{theorem}[Vitali]
Let $X$ be a measure space with finite measure $\mu,$ and let $h_k:
X\mapsto \mathbf C$ be a sequence of functions that is uniformly
integrable, i.e. such that for every $\varepsilon>0$ there exists
$\delta>0,$ independent of $k,$ satisfying $$ \mu(E) <\delta
\implies \int_E |h_k|\, d\mu <\varepsilon. \eqno(\dag)$$ Now: if\/
$\lim_{k\to\infty}h_k(x)=h(x)$ a.e., then
$$\lim_{k\to\infty}\int_X h_k\, d\mu=\int_X h\,d\mu.\eqno(\ddag)$$
In particular, if $$\sup_k\int_X |h_k|^p\,d\mu<\infty,\quad
\text{for some $p>1$},$$ then  $(\dag)$ and $(\ddag)$ hold.
\end{theorem}

\begin{proof}[Proof of Lemma~\ref{le1}]
For fixed $z\in {\mathbb B}$ let us make the change of variable
$w=\varphi_z(\omega),\ \omega\in {\mathbb B}$ in the integral which
represent $F_\zeta(z)$. In previous section we obtained the next
relation for pull-back measure
$$dv_\alpha(\varphi_z(\omega))=\frac{\left(1-|z|^2\right)^{n+1+\alpha}}{\left|1-\left<z,\omega\right>\right|^{2n+2+2\alpha}} dv_\alpha(\omega).$$
By using this result and \eqref{Relation2} we find
\[\begin{split}
\frac{F_\zeta(z)}{1+n+\alpha} &=\int_{\mathbb B} \frac
{(1-|z|^2)\left|\left<w,\overline{\zeta}\right>\right|}
{\left|1-\left<z,w\right>\right|^{n+2+\alpha}} dv_\alpha(w)
\\&=\int_{\mathbb B}
\frac
{(1-|z|^2)\left|\left<\varphi_z(\omega),\overline{\zeta}\right>\right|}
{\left|1-\left<z,\varphi_z(\omega)\right>\right|^{n+2+\alpha}}
\frac{\left(1-|z|^2\right)^{n+1+\alpha}}{\left|1-\left<z,\omega\right>\right|^{2n+2+2\alpha}}
dv_\alpha(\omega)
\\&=\int_{\mathbb B}
\frac
{(1-|z|^2)^{n+2+\alpha}\left|\left<\varphi_z(\omega),\overline{\zeta}\right>\right|}
{\left|1-\left<z,\varphi_z(\omega)\right>\right|^{n+2+\alpha}\left|1-\left<z,\omega\right>\right|^{2n+2+2\alpha}}
dv_\alpha(\omega)
\\&=\int_{\mathbb B}
\frac
{\left(|1-\left<z,\omega\right>|^{n+2+\alpha}|1-\left<z,\varphi_z(\omega)\right>|\right)^{n+2+\alpha}\left|\left<\varphi_z(\omega),\overline{\zeta}\right>\right|}
{\left|1-\left<z,\varphi_z(\omega)\right>\right|^{n+2+\alpha}\left|1-\left<z,\omega\right>\right|^{2n+2+2\alpha}}
dv_\alpha(\omega)
\\&=\int_{\mathbb B}
\frac {\left|\left<\varphi_z(\omega),\overline{\zeta}\right>\right|}
{\left|1-\left<z,\omega\right>\right|^{n+\alpha}} dv_\alpha(\omega).
\end{split}\]
Therefore \begin{equation}\label{popo}\frac{F_\zeta(z)}{1+n+\alpha}
=\int_{\mathbb B} \frac
{\left|\left<\varphi_z(\omega),\overline{\zeta}\right>\right|}
{\left|1-\left<z,\omega\right>\right|^{n+\alpha}}
dv_\alpha(\omega).\end{equation} From the last representation of
$F_\zeta(z)$ it follows
$$\frac{F_\zeta(z)}{1+n+\alpha}=\int_{\mathbb B}
\frac {\left|\left<\varphi_z(\omega),\overline{\zeta}\right>\right|}
{\left|1-\left<z,\omega\right>\right|^{n+\alpha}} dv_\alpha(\omega)
\le c_\alpha \int_{\mathbb B} \frac {(1-|\omega|^2)^\alpha
dv(\omega)}{\left|1-\left<z,\omega\right>\right|^{n+\alpha}}=
c_\alpha J_{c,t}(z),$$ where we set $c=-1$ and $t=\alpha$. Then
$\lambda_1=\frac{n+\alpha}2$, and $c_\alpha=\binom{n+\alpha}{n}$ as
in \eqref{cal}. For $z\in {\mathbb B},\ z\ne 0$ we have
$$J_{c,t}(z)\le J_{c,t}(z/|z|)=\frac{\Gamma(1 + n) \Gamma(1 + \alpha)}{\Gamma^2(
  (2 + n + \alpha)/2)}.$$ Thus
\begin{equation}\label{dri}\frac{F_\zeta(z)}{1+n+\alpha}\le \frac{c_\alpha \Gamma(1 + n) \Gamma(1 + \alpha)}{\Gamma^2(
 (2 + n + \alpha)/2)}=\frac{C_{\alpha,n}}{1+n+\alpha}
\end{equation}
what is the statement of lemma.
\end{proof}

\begin{proof}[Proof of Lemma~\ref{le2}]
Take $\zeta=e_1$ and $z=z_k=\frac{k}{k+1}\zeta$. Define
$$g_k(w)=\frac{w_1}{|w_1|}\frac{|1-\left<z,w\right>|^{n+2+\alpha}}{(1-\left<w,z\right>)^{n+2+\alpha}},\quad w\in {\mathbb B},\  w_1\ne 0.$$
Then $g_k\in L^\infty({\mathbb B})$ and $\|g_k\|_\infty=1$. Further
from \eqref{popo} and \eqref{nabla} we obtain
\begin{equation}\label{reso}
\begin{split}(1-|z_k|^2)|\nabla(
P_\alpha g_k)(z_k)|&\ge(1-|z_k|^2)|\left<\nabla( P_\alpha
g_k)(z_k),\zeta\right>|\\&=(1-|z_k|^2)\left|\int_{\mathbb B}
\left<\nabla_z
\mathcal{K}_\alpha(z,w)g_k(w),\zeta\right>dv_\alpha(w)\right|\\&=
(1+n+\alpha)\int_{\mathbb B}
\frac{(1-|z_k|^2)|w_1|}{\left|1-\left<z_k,w\right>\right|^{n+2+\alpha}}dv_\alpha(w)\\&=(1+n+\alpha)\int_{\mathbb
B} \frac
{\left|\left<\varphi_{z_k}(\omega),\overline{\zeta}\right>\right|}
{\left|1-\left<z_k,\omega\right>\right|^{n+\alpha}}
dv_\alpha(\omega):=G_k
\end{split}
\end{equation}
For $p=\frac{n +\alpha+1/2}{n+\alpha}$  ($p>1$), according to
Proposition~\ref{rudin} (take $c=-1/2$ and $t=\alpha$)
$$\sup_{k}\int_{\mathbb B} \left(\frac
{\left|\left<\varphi_{z_k}(\omega),\overline{\zeta}\right>\right|}
{\left|1-\left<z_k,\omega\right>\right|^{n+\alpha}}
\right)^pdv_\alpha(\omega)<\infty $$ (notice that
$\left|\left<\varphi_{z_k}(\omega),\overline{\zeta}\right>\right|\le
1$). Therefore by Vitali theorem
\[\begin{split}
\lim_{k\to \infty}G_k &=({1+n+\alpha})\lim_{k\to
\infty}\int_{\mathbb B} \frac
{\left|\left<\varphi_{z_k}(\omega),\overline{\zeta}\right>\right|}
{\left|1-\left<z_k,\omega\right>\right|^{n+\alpha}}
dv_\alpha(\omega)\\&=(1+n+\alpha)\int_{\mathbb B} \lim_{k\to
\infty}\frac
{\left|\left<\varphi_{z_k}(\omega),\overline{\zeta}\right>\right|}
{\left|1-\left<z_k,\omega\right>\right|^{n+\alpha}}
dv_\alpha(\omega).
\end{split}\]
For fixed $w\in {\mathbb B}$ we have
$$\lim_{k\to \infty}\frac
{\left|\left<\varphi_{z_k}(\omega),\overline{\zeta}\right>\right|}
{\left|1-\left<z_k,\omega\right>\right|^{n+\alpha}}=\frac{\left|\left<\zeta,\overline{\zeta}\right>\right|}
{\left|1-\left<\zeta,\omega\right>\right|^{n+\alpha}}=\frac{1}{{\left|1-\left<\zeta,\omega\right>\right|^{n+\alpha}}}.$$
Therefore by using Proposition~\ref{rudin} again we have
\[\begin{split}\lim_{k\to \infty}(1-|z_k|^2)|\nabla( P_\alpha
g)(z_k)|&=c_\alpha(1+n+\alpha)J_{-1,\alpha}(e_1)\\&=c_\alpha
\frac{(1+n+\alpha)\Gamma(1 + n) \Gamma(1 + \alpha)}{\Gamma^2(
 1/2 (2 + n + \alpha))}\\&=\frac{\Gamma(2 + n+\alpha)
}{\Gamma^2((2+n+\alpha)/2)}=C_{\alpha,n}.\end{split}\]
\end{proof}
\section{Proof of Theorem~\ref{Mait}}
Let
$$\tilde C_{\alpha,n}:=\|P_\alpha\|_{\tilde\beta}=\sup \{ |\tilde \nabla_z(P_\alpha g)(z)|  : |z|<1,\ \|g\|_\infty\le 1
\}.$$ We first prove \eqref{esti1}. It follows from the following
two lemmas.
\begin{lemma}
For $\alpha>-1$ and $\ell$ defined in \eqref{ell} we have $$\tilde
C_{\alpha,n}\le\max_{0\le t\le \pi/2}\ell(t).$$
\end{lemma}
\begin{proof} Let $f=P_\alpha g$.
We have
\[\begin{split}
(f\circ \varphi_a)(z)&=(P_\alpha g \circ \varphi_a)(z)
\\&=\int_{\mathbb B} K_\alpha(\varphi_a(z),w) g(w)dv_\alpha(w)
\\&=\int_{\mathbb B} K_\alpha(\varphi_a(z),\varphi_a(w)) g(\varphi_a(w))dv_\alpha(\varphi_a(w))
\end{split}\]
Since
$$1-\left<\varphi_a(z),\varphi_a(w)\right>=\frac{(1-\left<a,a\right>)(1-\left<z,w\right>)} {(1-\left<z,a\right>)(1-\left<a,w\right>)}$$
for \begin{equation}\label{thet}\theta=n+1+\alpha\end{equation} we
have
$$f\circ
\phi_a(z)=\frac{(1-\left<z,a\right>|)^\theta}{(1-|z|^2)^\theta}\int_{\mathbb
B}\frac{(1-\left<a,w\right>|)^\theta}{(1-\left<z,w\right>)^\theta}g\circ\varphi_a(w)
dv_\alpha(\varphi_a(w))$$ Differentiating in $z$ at 0 by using the
product rule we have
\begin{equation*}\label{poce}
\tilde\nabla f(a)=\theta\int_{\mathbb B}\frac{(\bar w-\bar
a)(1-\left<a,w\right>)^\theta}{(1-|a|^2)^\theta}g\circ
\varphi_a(w)dv_\alpha(\varphi_a(w)),
\end{equation*} where $$dv_\alpha(\varphi_a(w))=\left(\frac{\left(1-|a|^2\right)}{\left|1-\left<w,a\right>\right|^2}\right)^{n+1+\alpha} dv_\alpha(w).$$
Thus \begin{equation}\label{poca} \tilde\nabla
f(a)=\theta\int_{\mathbb B}\frac{(\bar w-\bar
a)(1-\left<a,w\right>)^\theta}{(1-|a|^2)^{\theta-n-1-\alpha}}\frac{g\circ
\varphi_a(w)}{|1-\left<w,a\right>|^{2(n+1+\alpha)} }dv_\alpha(w),
\end{equation}
and consequently for
\begin{equation}\label{thetap}\theta'=\binom{n+\alpha}{n}\theta\end{equation}
\begin{equation*}\label{poc}\begin{split} |\tilde\nabla
f(a)|&=\theta\sup_\zeta\left|\int_{\mathbb B}\left<\frac{(\bar
w-\bar
a)(1-\left<a,w\right>)^\theta}{(1-|a|^2)^{\theta-n-1-\alpha}}\frac{g\circ
\varphi_a(w)}{|1-\left<w,a\right>|^{2(n+1+\alpha)} },\zeta\right>
dv_\alpha(w)\right|\\&\le\theta'\sup_\zeta\int_{\mathbb
B}\left|\left<\frac{(\bar w-\bar
a)}{{|1-\left<w,a\right>|^{n+1+\alpha}}}
,\zeta\right>\right||{g\circ \varphi_a(w)}|(1-|w|^2)^\alpha
dv(w)\\&\le \theta'\|g\|_\infty\sup_\zeta\int_{\mathbb
B}\left|\left<{\bar w-\bar a} ,\zeta\right>\right|
\frac{(1-|w|^2)^\alpha}{{|1-\left<w,a\right>|^{n+1+\alpha}}} dv(w)
.\end{split}
\end{equation*}
Let $$L(a)=\sup_\zeta\int_{\mathbb B}\left|\left<{\bar w-\bar a}
,\zeta\right>\right|
\frac{(1-|w|^2)^\alpha}{{|1-\left<w,a\right>|^{n+1+\alpha}}} dv(w)$$
and define
$$L=\sup_{a\in\mathbb B}L(a).$$ Then $$L=\sup_{a\in\mathbb B}\sup_\zeta\int_{\mathbb B}\left|S_{\zeta,w}(a)\right|
dv_\alpha(w),$$ where
$$S_{\zeta,w}(a)=\frac{\left<w-a,\zeta\right>}{{(1-\left<w,a\right>)}^{n+1+\alpha}}.$$
Observe that $S_{\zeta,w}(a)$ is a subharmonic function in $a$. It
follows that $a\to L(a)$ is subharmonic and its maximum is achieved
on the boundary of the unit ball. Therefore there exist
$a_0,\zeta_0\in \mathbb S$ such that
$$L=\int_{\mathbb B}\left|\left<{ w-a_0}
,\zeta_0\right>\right|
\frac{(1-|w|^2)^\alpha}{{|1-\left<w,a_0\right>|^{n+1+\alpha}}}
dv(w).$$ Let $U$ be an unitary transformation of $\mathbb{C}^n$ onto
itself such that $Ua_0=e_1$ and $U\zeta_0=\cos t e_1+\sin t e_2$ for
some $t\in[0,\pi]$ (Here $t=\arg(a_0,\zeta_0)$). Take the
substitution $w=U\omega$. Then we obtain
\[\begin{split}L&=\int_{\mathbb B}\left|\left<{ U\omega-Ue_1} ,\zeta_0\right>\right|
\frac{(1-|U\omega|^2)^\alpha}{{|1-\left<U\omega,a_0\right>|^{n+1+\alpha}}}
dv(U\omega)\\&=\int_{\mathbb B}\left|\left<{ \omega-e_1}
,U\zeta_0\right>\right|
\frac{(1-|\omega|^2)^\alpha}{{|1-\left<\omega,Ua_0\right>|^{n+1+\alpha}}}
dv(\omega)\\&=\int_{\mathbb B}\left|\left<{ \omega-e_1} ,\cos t
e_1+\sin t e_2\right>\right|
\frac{(1-|\omega|^2)^\alpha}{{|1-\left<\omega,e_1\right>|^{n+1+\alpha}}}
dv(\omega)\\&=\int_{\mathbb B}\frac{|(1-w_1)\cos t+w_2\sin
t|}{|w_1-1|^{n+1+\alpha}}dv_\alpha(w).\end{split}\]

\end{proof}
\begin{lemma}\label{22} Let $\ell$ be defined as in \eqref{ell}. Then
$$\tilde C_{\alpha,n}\ge \ell(\pi/2).$$
\end{lemma}
\begin{proof} Let $\zeta=e_2$, $a=\epsilon_k e_1$, where $\epsilon_n =k/(k+1)$. Then
$$|\tilde \nabla f(a)\zeta|=\theta'\left|\int_{\mathbb B}{\bar w_2(1-\epsilon_k w_1)^{n+1+\alpha}}{g\circ
\varphi_a(w)}\frac{(1-|w|^2)^\alpha}{{|1-\epsilon_k
w_1|^{2(n+1+\alpha)} }} dv(w)\right|.$$ Define $g_k$ such that
$$\bar w_2(1-\epsilon_k w_1)^{n+1+\alpha}{g_k\circ \varphi_a(w)}=\left|\bar w_2(1-\epsilon_k w_1)^{n+1+\alpha}\right|$$ and let
$f_k=P_\alpha(g_k)$. Then we have $$|\tilde \nabla
f_k(a)\zeta|=\theta'\int_{\mathbb B}{
|w_2|}\frac{(1-|w|^2)^\alpha}{{|1-\epsilon_k w_1|^{n+1+\alpha} }}
dv(w).$$ Thus $$\tilde C_{\alpha,n}\ge \sup_{k,\zeta,a}|\tilde
\nabla f_k(a)\zeta|\ge \ell(\pi/2).$$
\end{proof}
In order to prove \eqref{elles} and \eqref{esti2} we prove the
following lemma (which is an extension of a corresponding result of
Bungart, Folland and Fefferman, cf. Rudin book
\cite[Proposition~1.4.9]{Rudin-Book-Ball}).
\begin{lemma}\label{lepo} For a multi-index $\eta=(\eta_1,\dots,\eta_n)\in\mathbb{N}^n_0$
we have \begin{equation}\label{generud}\int_S
|\zeta^\eta|d\sigma(\zeta)=\frac{(n-1)!\prod_{i=1}^n\Gamma[1+\frac{\eta_i}{2}]}{\Gamma[n+\frac{|\eta|}{2}]}\end{equation}
and \begin{equation}\label{generPO}\int_{\mathbb B}
|z^\eta|dv_\alpha(z)=\frac{\Gamma[1+\alpha+n]}{
\Gamma[1+\alpha+|\eta|/2+n]}{\prod_{i=1}^n\Gamma[1+\frac{\eta_i}{2}]}.\end{equation}
Here $w^\eta:=\prod_{k=1}^n w_k^{\eta_k}$ and
$|\eta|=\sum_{i=1}^n\eta_i$.
\end{lemma}
Notice that the following proof works as well assuming that
$\eta=(\eta_1,\dots,\eta_n)$, where $\eta_j>-1$, $j=1,\dots,n$ are
arbitrary real numbers.
\begin{proof}[Proof of Lemma~\ref{lepo}]
The proof of \eqref{generud} goes along the proof of the similar
statement in Rudin's book where it is proved the same statement  for
$\eta=2\chi$, where $\chi$ is a multi-index. Here are details for
the sake of completeness. Let
$$I=\int_{\mathbb C^n}|z^\eta|\exp(-|z|^2)dV(z),$$
where $dV$ is the Lebesgue measure in $\mathbb C^n$. The expression
under integral is
$$\prod_{j=1}^n|z_j|^{\eta_j}\exp(-|z|^2).$$
By Fubini's theorem
$$I=\prod_{j=1}^n\int_\mathbb C|\lambda|^{\eta_j}\exp(-|\lambda|^2)dV(\lambda).$$
One can easy compute the next for $m,\ \Re(m)>-1$
$$\int_\mathbb C|\lambda|^m\exp(-|\lambda|^2)dV(\lambda).$$
Namely, by using polar coordinates $\lambda=r\zeta,\ \zeta\in\mathbb
T$, we obtain
\[\begin{split}
\int_\mathbb
C|\lambda|^m\exp(-|\lambda|^2)dv(\lambda)&=2\int_0^\infty
rdr\int_\mathbb T |r\zeta|^m\exp(-|r\zeta|^2)d\sigma(\zeta)
\\&=2\int_0^\infty rdr\int_\mathbb T r^m\exp(-r^2)d\sigma(\zeta)
\\&=2\int_0^\infty rdr\int_\mathbb T r^m\exp(-r^2)d\sigma(\zeta)
\\&=2\int_0^\infty r^{m+1}\exp(-r^2)dr
\\&=\int_0^\infty (r^2)^{m/2}\exp(-r^2)d(r^2)
\\&=\int_0^\infty t^{m/2}\exp(-t)dt
\\&=\Gamma(1 +m/2 ).
\end{split}\]
Thus
$$\int_\mathbb C|\lambda|^m\exp(-|\lambda|^2)dV(\lambda)=\pi \Gamma(1 +m/2 ).$$
and it follows
$$I= \pi ^n \prod_{j=1}^n\Gamma(1 +\eta_j/2 ).$$
On the other hand, applying polar coordinates in $I$, we obtain
($\omega_{2n}$  is  volume measure of unit ball)
$$ I/\omega_{2n} =2n \int_0^\infty r^{|\eta|+2n-1}\exp(-r^2)dr\int_S|\zeta^\eta| d\sigma(\zeta).$$
Thus
$$\int_S|\zeta^\eta| d\sigma(\zeta)= I/(\omega_{2n} \cdot 2n \int_0^\infty r^{|\eta|+2n-1}\exp(-r^2)dr).$$
Since $\omega_{2n}=\pi^n/n!$ and $2\int_0^\infty
r^{|\eta|+2n-1}\exp(-r^2)dr=\Gamma(n+|\eta|/2)$ it follows
$$\int_S|\zeta^\eta| d\sigma(\zeta)=(n-1)!\frac{ \prod_{j=1}^n\Gamma(1 +\eta_j/2 )}{\Gamma(n+|\eta|/2)} .$$

Let us prove now \eqref{generPO}. For a mapping $f\in L^1(\mathbb
B)$ we have
\[\begin{split}\int_{\mathbb B}f(x) dv_\alpha(x)&=\binom{n+\alpha}{n}\int_{\mathbb
B}(1-|x|^2)^\alpha f(x)
dv(x)\\&=2n\binom{n+\alpha}{n}\int_0^1r^{2n-1}(1-r^2)^\alpha
\int_{\mathbb S} f(r\eta)d\sigma(\eta)dr .\end{split}\]

For $f(z)=|z|^\eta$ we have $$\int_{\mathbb S}
f(r\zeta)d\sigma(\zeta)=r^{|\eta|}\frac{(n-1)!\prod_{i=1}^n\Gamma[1+\frac{\eta_i}{2}]}{\Gamma[n+\frac{|\eta|}{2}]}.$$

Further we have
$$2n\binom{n+\alpha}{n}\int_0^1r^{2n+|\eta|-1}(1-r^2)^\alpha dr =\frac{\Gamma[1+\alpha+n]\Gamma[|\eta|/2+n]}{\Gamma(n) \Gamma[1+\alpha+|\eta|/2+n]}.$$
This finishes the proof of the lemma.
\end{proof}
The relation \eqref{elles}  and the left-hand inequality in
\eqref{esti2} follows from the following lemma (in view of
\eqref{elles}).
\begin{lemma}\label{ika}
Let $\ell(t)$ be defined as in \eqref{ell}. Then
$\ell(0)=C_{\alpha,n}$ and $\ell(\pi/2)=\frac{\pi}{2}C_{\alpha,n}$.
\end{lemma}
\begin{proof}[Proof of Lemma~\ref{ika}] The relation $\ell(0)=C_{\alpha,n}$ follows at once. Prove the second relation.

Observe first that  for $l\neq k$
$$\int_{\mathbb B} w_{1}^{l}\bar{w}_1^{k}|w_2|dv(w)=0.$$ By choosing
$\eta(k)=(2k,1,0,\dots,0)$ we obtain
\[\begin{split}J&=\int_{\mathbb B}\frac{|w_2|}{|1-w_1|^{n+1+\alpha}}dv_\alpha(w)\\&=\int_{\mathbb B}\frac{|w_2|}{|(1-w_1)^{(n+1+\alpha)/2}|^2}dv_\alpha(w)\\&=\sum_{k=0}^\infty
\binom{-(n+1+\alpha)/2}{k}^2\int_{\mathbb
B}|w_1|^{2k}|w_2|dv_\alpha(w),\ \ \ w=(w_1,\dots,w_n).\end{split}\]
From \eqref{generPO} we find that
$$\label{gener} \int_{\mathbb
B}|w_1|^{2k}|w_2|dv_\alpha(w)=\int_{\mathbb
B}|z^{\eta(k)}|dv_\alpha(z)=\frac{\Gamma[1+\alpha+n]}{
\Gamma[1+\alpha+|\eta|/2+n]}{\Gamma[3/2]\Gamma[1+k]}.$$ Therefore
\[\begin{split}J&=\Gamma[3/2]\Gamma[1+\alpha+n]\sum_{k=0}^\infty
\binom{-(n+1+\alpha)/2}{k}^2\frac{{k!}}{
\Gamma[1+\alpha+k+n+1/2]}\\&=\frac{\Gamma[3/2]\Gamma[1+\alpha+n]}{\Gamma[\alpha+n+3/2]}\sum_{k=0}^\infty
\frac{(((n+1+\alpha)/2)_k)^2}{ (\alpha+n+3/2)_k
k!}\\&=\frac{\Gamma[3/2]\Gamma[1+\alpha+n]}{\Gamma[\alpha+n+3/2]}F[(n+1+\alpha)/2,(n+1+\alpha)/2,(\alpha+n+3/2),1]\\&=\frac{\pi
\Gamma[1 + a + n]}{2 \Gamma[(2 + a + n)/2]^2}.\end{split}\] The last
equality is derived with help of Gauss theorem i.e. of the relation
\eqref{fora}. Hence
\[\begin{split}\ell(\pi/2)&=(n+1+\alpha)\int_{\mathbb B}\frac{|w_2|}{|1-w_1|^{n+1}}dv_\alpha(w)\\&=\frac{\pi
(n+1+\alpha)\Gamma[1 + a + n]}{2 \Gamma[(2 + a + n)/2]^2}=\frac{\pi
\Gamma[2 + a + n]}{2 \Gamma[(2 + a +
n)/2]^2}\\&\left(>\frac{\Gamma(2 + n+\alpha)
}{\Gamma^2((2+n+\alpha)/2)}=\|P\|_\beta\right).\end{split}\]
\end{proof}
To finish the the proof of Theorem~\ref{Mait} we need to prove the
right inequality in \eqref{esti2}. It follows from this simple
observation
$$\tilde C_{\alpha,n}\le |\sin t| \ell(0)+|\cos t|\ell(\pi/2)\le \sqrt{\ell(0)^2+\ell(\pi/2)^2}.$$

\begin{remark}
If $g\in C(\overline{\mathbb{B}})$ and $f=P_\alpha[g]$, then it
follows from \eqref{poca} and Vitali theorem that there exist a
mapping $\Phi:\mathbb{S}\to \mathbb{C}^n$ such that for $\zeta\in
\mathbb{S}$
$$\lim_{a\to \zeta}\tilde\nabla f(a)=g(\zeta)\Phi(\zeta).$$ But if $g$ is a
polinom, then we know that $\lim_{a\to t}\tilde\nabla f=0$, implying
that $$\Phi(\zeta)=\int_{\mathbb B}\frac{(\bar w-\bar
\zeta)(1-\left<\zeta,w\right>)^{n+1+\alpha}}{|1-\left<w,\zeta\right>|^{2(n+1+\alpha)}
}dv_\alpha(w)=0, \ \text{for}\ \  \zeta\in \mathbb{S}.$$
If by
$\mathcal{B}_0$ we denote the little Bloch space, i.e. the space of
holomorphic mappings $f$ defined on the unit ball such that
$$\lim_{|z|\to 1}|\tilde \nabla f(z)|=0, $$ and consider the Bergman
projection $$P_\alpha: C(\overline{\mathbb{B}})\to \mathcal{B},$$
then by the previous consideration we obtain
$$P_\alpha(C(\overline{\mathbb{B}}))\subset \mathcal{B}_0\subset
\mathcal{B}.$$ It follows that $$\| P_\alpha:
C(\overline{\mathbb{B}})\to \mathcal{B}_0\|\le
\tilde{C}_{\alpha,n}$$ w.r.t. invariant $\tilde \beta$ Bloch
semi-norm. Moreover, since the extremal sequence (see the proof of
Lemma~\ref{22}) is consisted of continuous functions $g_k$, it
follows that
$$\| P_\alpha: C(\overline{\mathbb{B}})\to \mathcal{B}_0\|=
\tilde{C}_{\alpha,n}.$$ The same can be repeated for the standard
$\beta$ Bloch semi-norm.
\end{remark}

\subsection{Appendix}
Regarding Conjecture~\ref{cok} we offer the following observation.

Put $\alpha=0$ and $n=2$.  Let
$$I(t)=\frac{\ell(t)}{(n+1+\alpha)}=\int_{\mathbb{B}}\frac{|(1-w_1)\cos t+w_2\sin
t|}{|w_1-1|^{3}}dv(w).$$ Take the substitution $\varphi$:
$a_1=w_2/(1-w_1)$, $a_2=w_1$ on the integral. We obtain
$$\varphi^{-1}(\mathbb{B})=B'=\{(a_1,a_2): |a_1|^2\le
\frac{1-|a_2|^2}{|1-a_2|^2}\},$$ and ${|w_1-1|^{-2}}dv(w)=dv(a)$ and
$$I(t)=\int_{B'}|\cos t + a_1\sin t|dv(a).$$
Let $a_2=pe^{is}$, $0<p<1$,
$$a_1=R e^{i \sigma}$$ $$0\le R\le R_0=\frac{1-p^2}{1+p^2-2p \cos s}.$$
Then $$\int_{B'}|\cos t + a_1\sin
t|dv(a)=\int_0^{2\pi}ds\int_0^{R_0} dR\int_0^1 dp \int_0^{2\pi}|\cos
t +a_1\sin t| J d\sigma,$$ where $J=pR$. Define
$$h(t)=\int_0^{2\pi}|\cos t +a_1\sin t|
d\sigma=\int_0^{2\pi}\sqrt{\cos^2 t + R \sin t (2 \cos\sigma \cos t
+ R \sin t)}d\sigma.$$ Before we go further remark the following,
$h(t)$ is indeed the circumference of the ellipse $E[a,b]$ with the
axis $a=\cos t + R \sin t$ and $b=|\cos t - R \sin t|$. It can be
expressed by the formula
$$h(t)=4aE[\epsilon^2]$$ where $E$ is the elliptic function of the
second kind and
$$\epsilon=\sqrt{1-b^2/a^2}=\sqrt{\frac{2R \sin(2t)}{(\cos t  + R
\sin t)^2}}$$ is the eccentricity of the ellipse.

Therefore $$h'(t)=\csc(2 t) K[\epsilon^2] (-\cos t+R \sin t)+\cot(2
t) E[\epsilon^2] (\cos t+R \sin t),$$ where $K$ is the elliptic
function of the first kind. By using the asymptotic formulas
$$K[\epsilon^2]\approx \frac{\pi}{2}+\frac{\pi}{8}{\epsilon^2}$$ and
$$E[\epsilon^2]\approx \frac{\pi}{2}-\frac{\pi}{8}\epsilon^2$$ we
obtain  that
$$h'(0):=\lim_{t\to 0}h'(t)=h'(\pi/2):=\lim_{t\to \pi/2}h'(t)=0.$$ This means that
$I'(0)=I'(\pi/2)=0$, or what is the same $0$ and $\pi/2$ are
stationary points of the function $I$.






\end{document}